\documentclass[12pt]{article}
\usepackage{epic,eepic,epsf,epsfig}
\usepackage{amsmath}
\usepackage{amssymb}
\usepackage{amsbsy}

\usepackage{amsfonts}%,srcltx}
\newtheorem{theorem}{Theorem}[section]%
\newtheorem{lemma}[theorem]{Lemma}%
\newtheorem{cor}[theorem]{Corollary}%

\def\mz{{\mathbb Z}}

\def\f{\noindent}

\def\Syl{\hbox{\rm Syl}}

\newcommand{\qed}{\mbox{\raisebox{0.7ex}{\fbox{}}} \vspace{4truemm}}
\def\demo{\f {\bf Proof.}\hskip10pt}

\begin{document}

\baselineskip 16pt

\title{ \vspace{-1.2cm}
Finite groups in which some particular non-nilpotent maximal invariant subgroups have indices a prime-power
\thanks{\scriptsize This research was supported in part by Shandong Provincial Natural Science Foundation, China (ZR2017MA022)
and NSFC (11761079).
\newline
 \hspace*{0.5cm} \scriptsize $^{\ast\ast}$Corresponding
  author.
\newline
       \hspace*{0.5cm} \scriptsize{E-mail addresses:}
       shijt2005@163.com\,(J. Shi),\,shanmj2023@s.ytu.edu.cn\,(M. Shan),\,xufj2023@s.ytu.edu.cn\,(F. Xu).}}

\author{Jiangtao Shi\,$^{\ast\ast}$,\,Mengjiao Shan,\,Fanjie Xu\\
\\
{\small School of Mathematics and Information Sciences, Yantai University, Yantai 264005, China}}

\date{ }

\maketitle \vspace{-.8cm}

\begin{abstract}
Let $A$ and $G$ be finite groups such that $A$ acts coprimely on $G$ by automorphisms, assume that $G$ has a
maximal $A$-invariant subgroup $M$ that is a direct product of some isomorphic simple groups,
we prove that if $G$ has a non-trivial $A$-invariant normal subgroup $N$ such that
$N\leq M$ and every non-nilpotent maximal $A$-invariant subgroup $K$ of $G$ not containing
$N$ has index a prime-power and the projective special linear group $PSL_2(7)$ is not a composition factor of $G$, then $G$ is solvable.

\medskip
\f {\bf Keywords:} non-nilpotent maximal invariant subgroup; invariant normal subgroup; prime-power index; solvable\\
{\bf MSC(2020):} 20D10
\end{abstract}

\section{Introduction}

All groups are assumed to be finite. It is known that a theorem of Hall {\rm\cite[Theorem 9.4.1]{rob}} indicated that a group $G$
is solvable if every maximal subgroup of $G$ has index a prime or the square of a prime. As a generalization  of Hall's theorem,
Shi, Shi and Zhang {\rm\cite[Theorem 13]{SZ}} gave the result:

\begin{theorem} {\rm\cite[Theorem 13]{SZ}}\ \  Let $M$ be a maximal subgroup of $G$ such that $M$ is isomorphic to a direct
product of some isomorphic simple groups. Assume that there exists a subgroup $1\neq N\leq M$ such that $N\unlhd G$ and for every maximal
subgroup $K$ of $G$ that does not contain $N$ we always have that $|G:K|$ is a prime. Then $G$ is solvable.
\end{theorem}

Moreover, Shi, Kutnar and Zhang {\rm\cite[Theorem 3.6]{shi2018}} proved the result:

\begin{theorem} {\rm\cite[Theorem 3.6]{shi2018}}\ \ Let $G$ be a group with a maximal subgroup $M$ being a characteristic simple group.
Suppose that $G$ has a non-trivial normal subgroup $N$ such that $N\leq M$ and that the index of every maximal subgroup $K$ of $G$ not containing
$N$ is either a prime or  a square of a prime. Then $G$ is solvable.
\end{theorem}

As an extension of {\rm\cite[Theorem 3.6]{shi2018}}, combine the coprime action of groups, Shi {\rm\cite[Theorem 1.6]{shi2023}} obtained the following result.

\begin{theorem} {\rm\cite[Theorem 1.6]{shi2023}}\ \ Let $A$ and $G$ be groups such that $A$ acts coprimely on $G$ by automorphisms, assume that $G$ has a
maximal $A$-invariant subgroup $M$ that is a direct product of some isomorphic simple groups. If $G$ has a non-trivial $A$-invariant normal subgroup $N$ such that
$N\leq M$ and every non-nilpotent maximal $A$-invariant subgroup $K$ of $G$ not containing $N$ has index a prime or the square of a prime, then $G$ is solvable.
\end{theorem}

Consider maximal subgroups of a group having indices a prime-power, Guralnick {\rm\cite[Corollary 3]{gur}} proved that a group $G$ is solvable or $G/O_{\infty}G\cong PSL_2(7)$ if all maximal
subgroups of $G$ have indices a prime-power, where $O_{\infty}G$ is the largest normal solvable
subgroup of $G$. As a generalization of Guralnick's result, Guo {\rm\cite[Theorem 1]{guo}} indicated that $G/S(G)\cong 1$ or $PSL_2(7)$ if all non-nilpotent maximal subgroups of $G$ have indices a prime-power, where $S(G)$ is the largest normal solvable subgroup of $G$.

Combine the coprime action of groups, consider non-nilpotent maximal invariant subgroups of a group having indices a prime-power, Beltr$\rm\acute{a}$n and Shao {\rm\cite[Theorem D]{beltran}} gave the result:

\begin{theorem} {\rm\cite[Theorem D]{beltran}}\ \   Suppose that a group $A$ acts coprimely on a group $G$ and
let $p$ be a prime divisor of $|G|$. If the indices of all non-nilpotent maximal $A$-invariant
subgroups of $G$ are powers of $p$, then $G$ is solvable.
\end{theorem}

Recently, Shi and Liu {\rm\cite[Theorem 1.8]{shiliu2024}} provided the result:

\begin{theorem} {\rm\cite[Theorem 1.8]{shiliu2024}}\ \  Suppose that $A$ and $G$ are groups such that $A$ acts coprimely on $G$ by automorphisms. If every non-nilpotent maximal $A$-invariant
subgroup of $G$ has index a prime-power and $PSL_2(7)$ is not a composition factor of $G$, then $G$ is solvable.
\end{theorem}

In this paper, as a further extension of {\rm\cite[Theorem 1.6]{shi2023}} and a further generalization of {\rm\cite[Theorem D]{beltran}} and {\rm\cite[Theorem 1.8]{shiliu2024}}, we have the following result
whose proof is given in Section~\ref{s3}.

\begin{theorem} \ \ \label{th1} Let $A$ and $G$ be groups such that $A$ acts coprimely on $G$ by automorphisms, assume that $G$ has a maximal $A$-invariant subgroup $M$
that is a direct product of some isomorphic simple groups.
If $G$ has a non-trivial $A$-invariant normal subgroup $N$ such that $N\leq M$ and every non-nilpotent maximal $A$-invariant subgroup $K$ of $G$ not containing
$N$ has index a prime-power and $PSL_2(7)$ is not a composition factor of $G$, then $G$ is solvable.
\end{theorem}

By Case I in the proof of Theorem~\ref{th1}, we also obtain the result:

\begin{theorem} \ \ \label{th1} Let $A$ and $G$ be groups such that $A$ acts coprimely on $G$ by automorphisms, assume that $G$ has a maximal $A$-invariant subgroup $M$
that is a direct product of some isomorphic simple groups.
If $G$ has a non-trivial $A$-invariant normal subgroup $N$ such that $N\leq M$ and every maximal $A$-invariant subgroup $K$ of $G$ not containing
$N$ is nilpotent, then $G$ is solvable.
\end{theorem}

When $A=1$, the following corollary emerges from Theorem~\ref{th1}.

\begin{cor} \ \ \label{c1} Let $G$ be a group having a maximal subgroup $M$ that is a direct product of some isomorphic simple groups.
If $G$ has a non-trivial normal subgroup $N$ such that $N\leq M$ and every non-nilpotent maximal subgroup $K$ of $G$ not containing
$N$ has index a prime-power and $PSL_2(7)$ is not a composition factor of $G$, then $G$ is solvable.
\end{cor}

\section{Some Necessary Lemmas}

\begin{lemma} {\rm\cite[Theorem B]{beltran}} \ \ \label{l1} Let $G$ and $A$ be groups of coprime orders and assume that $A$ acts on $G$ by
automorphisms. If $G$ has a nilpotent maximal $A$-invariant subgroup of odd order, then $G$ is
solvable.
\end{lemma}

\begin{lemma} {\rm\cite[Theorem D]{beltran}} \ \label{l2}  Suppose that a group $A$ acts coprimely on a group $G$ and
let $p$ be a prime divisor of $|G|$. If the indices of all non-nilpotent maximal $A$-invariant
subgroups of $G$ are powers of $p$, then $G$ is solvable.
\end{lemma}

\begin{lemma} {\rm\cite[Lemma 2.2]{shi2023}}\ \ \label{l3} Suppose that $A$ and $G$ are groups such that $A$ acts coprimely on $G$ by automorphisms.
If $G$ has an abelian maximal $A$-invariant subgroup $L$, then $G$ is solvable.
\end{lemma}

\begin{lemma}  {\rm\cite[Lemma 3.3]{shi2018}}\ \ \label{l4}
Let $G$ be a non-abelian simple group. Suppose that
there exist $H<G$ and $K<G$ such that $|G:H|=p^m$ and $|G:K|=q^n$,
where $p$ and $q$ are primes, and $m$ and $n$ are positive integers.
Then one of the following holds:

$(1)$ $p^m=q^n$;

$(2)$ $G=PSL_2(7),\,H\cong S_4,\,K\cong\mz_7\rtimes\mz_3$.
\end{lemma}

\begin{lemma} {\rm\cite[Lemma 3.4]{shi2018}}\ \ \label{l5} Let $G=G_1\times
G_2\times\cdots\times G_s$, where $G_i\cong G_j$ is a non-abelian
simple group for every $i,j\in\{1,\cdots,s\}$. If there exists $K<G$
such that $|G:K|=p^n$, $p$ a prime and $n$   a positive
integer, then for every $i\in\{1,\cdots,s\}$ there exists $M_i<G_i$ such that $|G_i:M_i|=p^m$,
where $m$ is a positive integer.
\end{lemma}

\begin{lemma} {\rm\cite[Lemma 9]{wie}}\ \ \label{l6} Let $H$ be a nilpotent Hall-subgroup of a group $G$ that is not a Sylow subgroup
of $G$. If for each prime divisor $p$ of $|H|$, assume $P\in\Syl_p(H)$ we always have $N_G(P)=H$, then there exists a normal
subgroup $K$ of $G$ such that $G=KH$ and $K\cap H=1$.
\end{lemma}

\section{Proof of Theorem~\ref{th1}}\label{s3}

\demo By the hypothesis, let $N_0$ be a minimal $A$-invariant normal subgroup of $G$ such that $N_0\leq N$. Then for any non-nilpotent maximal $A$-invariant
subgroup $K$ of $G$, if $N_0\nleq K$, one has that $N\nleq K$ and $|G:K|$ is equal to a prime-power.

Assume that $G$ is non-solvable and let $G$ be a counterexample of minimal order. By Lemma~\ref{l3}, $M$ is non-abelian. It follows that $M$ is a direct
product of some isomorphic non-abelian simple groups. Moreover, $N_0$ is also a direct product of some isomorphic non-abelian simple groups by
{\rm\cite[Proposition 1.6.3]{kur}}. Assume $N_0=N_{01}\times N_{02}\times\cdots\times N_{0s}$, where $s\geq 1$ and $N_{01},\,N_{02}\,,\cdots,\,N_{0s}$
are isomorphic non-abelian simple groups.

Let $p$ be the largest prime divisor of $|N_0|$ and $P$ an $A$-invariant Sylow $p$-subgroup of $G$. Then $P_0=P\cap N_0$ is an $A$-invariant
Sylow $p$-subgroup of $N_0$. By the minimality and non-solvability of $N_0$, one has that $P_0$ is not normal in
$G$, which implies $N_G(P_0)<G$. Then $G=N_0N_G(P_0)$ by Frattini's argument. Let
$L$ be a maximal $A$-invariant subgroup of $G$ such that $N_G(P_0)\leq L$. One has $G=N_0L$, which implies $N_0\nleq L$.

For $L$, we divide our analyse into two cases.

{\bf Case I.} When $L$ is nilpotent.

First, claim that $L$ is a Hall-subgroup of $G$. Otherwise, if $L$ is not a Hall-subgroup of $G$, then $L$ has an $A$-invariant Sylow subgroup $T$ that
is not a Sylow subgroup of $G$. One has $L<N_G(T)$. Since $N_G(T)$ is also an $A$-invariant subgroup of $G$, $N_G(T)=G$ by the maximality of $L$. It follows
that $T$ is normal in $G$. Let $T_0$ be a minimal $A$-invariant normal subgroup of $G$ such that $T_0\leq T$. It is obvious that $T_0\neq N_0$. Consider
the quotient group $G/T_0$. It is easy to see that the hypothesis of the theorem holds for $G/T_0$. By the minimality of $G$, one has that $G/T_0$ is solvable,
which implies that $G$ is solvable, a contradiction. Therefore, $L$ is a Hall-subgroup of $G$.

Second, claim that $L$ is not a Sylow subgroup of $G$. Otherwise, assume that $L$ is a Sylow subgroup of $G$. Since $p$ is the largest prime
divisor of $|N_0|$ and $N_0$ is non-solvable, $p$ is an odd prime. It follows that $L=P_0$ is a nilpotent maximal $A$-invariant subgroup of $G$ of odd order
since $L\geq P_0$, which implies that $G$ is solvable by Lemma~\ref{l1}, a contradiction. Therefore, $L$ is not a Sylow subgroup of $G$.

Third, let $W_0$ be any $A$-invariant Sylow subgroup of $L$. Arguing as above $T_0$, $G$ has no solvable $A$-invariant normal subgroups. Since $N_G(W_0)$ is
also an $A$-invariant subgroup of $G$ and $L\leq N_G(W_0)<G$, one has $N_G(W_0)=L$ by the maximality of $L$.

Therefore, there exists a normal subgroup $H$ of $G$ such that $G=H\rtimes L$ by Lemma~\ref{l6}. Note that $L$ is a Hall-subgroup of $G$, which implies that $H$ is a normal
Hall-subgroup of $G$. It follows that $H$ is a characteristic subgroup of $G$
and then $H$ is an $A$-invariant normal subgroup of $G$. Moreover, $H$ is a minimal $A$-invariant normal subgroup of $G$ by the maximality of $L$. It is
easy to see that $H\neq N_0$ since $p\nmid|H|$ and $p\mid|N_0|$.
It follows that $H\cap N_0=1$. It is easy to see that $G\cong G/(H\cap N_0)\lesssim G/H\times G/N_0$, $G/H\cong L$ is nilpotent and $G/N_0=N_0L/N_0\cong
L/(N_0\cap L)$ is also nilpotent, one has that $G$ is nilpotent, a contradiction.

{\bf Case II.} When $L$ is non-nilpotent.

By the hypothesis, $|G:L|$ is equal to a prime-power. Since $|G:L|=|N_0L:L|=|N_0:N_0\cap L|$ and $N_0\cap L\geq N_0\cap N_G(P_0)=N_{N_0}(P_0)$,
one has $|G:L|=q^{\alpha}$, where $q\neq p$ is a prime divisor of $|N_0|$ and $\alpha\geq 1$. It follows that $N_0$ has a proper subgroup of index a $q$-power, which
implies that $N_{01}$ has a proper subgroup of index a $q$-power by Lemma~\ref{l5}.

Let $Q$ be an $A$-invariant Sylow $q$-subgroup of $G$. Then $Q_0=Q\cap N_0$ is an $A$-invariant Sylow $q$-subgroup of $N_0$. Arguing as above, $Q_0$ is not
normal in $G$ and $G=N_0N_G(Q_0)$. Let $R$ be a maximal $A$-invariant
subgroup of $G$ such that $N_G(Q_0)\leq R$. Then $G=N_0R$.

For $R$, we divide our arguments into two cases.

{\bf Case $(i)$.} When $R$ is nilpotent.

Arguing as above, one has that $R$ is a Hall-subgroup of $G$ and for any $A$-invariant Sylow subgroup $U$ of $R$ we have $N_G(U)=R$.
When $R$ is not a Sylow subgroup of $G$ or $R$ is a Sylow subgroup of $G$ of odd order, we can always get a contradiction arguing as in Case I. In the following, assume that $R$ is an
$A$-invariant Sylow 2-subgroup of $G$, which implies $q=2$ and $|G:L|=2^{\alpha}$. It follows that $|N_0:N_0\cap L|=2^{\alpha}$. By Lemma~\ref{l5},
$N_{01}$ has a proper subgroup of index an 2-power.

If all non-nilpotent maximal $A$-invariant subgroups of $G$ have indices an $2$-power. Then $G$ is solvable by Lemma~\ref{l2}, a contradiction.

Therefore, there exists a non-nilpotent maximal $A$-invariant subgroup $V$ of $G$ such that $|G:V|$ is not equal to an $2$-power. If $N_0\leq V$. Then $G=N_0R=VR$, which
implies $|G:V|=|VR:V|=|R:V\cap R|$ is an 2-power, a contradiction. If $N_0\nleq V$. Then $|G:V|$ is a prime-power by the hypothesis. Assume $|G:V|=k^{\beta}$,
where $k\neq 2$. It follows that $|N_0:N_0\cap V|=k^{\beta}$ and then $N_{01}$ has a proper subgroup of index a $k$-power by Lemma~\ref{l5}. Note that $k\neq 2$. One has
$N_{01}\cong PSL_2(7)$ by Lemma~\ref{l4}, which implies that $G$ has a composition factor that is isomorphic to $PSL_2(7)$, a contradiction.

{\bf Case $(ii)$.} When $R$ is non-nilpotent.

Since $N_0\nleq R$, $|G:R|=r^{\tau}$ by the hypothesis, where $r\neq q$ is a prime divisor of $|N_0|$ and $\tau\geq 1$. It
implies that $N_{01}$ has a proper subgroup of index a $r$-power by Lemma~\ref{l5}. Arguing as above, we can get $N_{01}\cong PSL_2(7)$, also a contradiction.

So the counterexample of minimal order does not exist and then $G$ is solvable.\hfill\qed

\end{document}